\documentclass[12pt]{article}%
\usepackage{amsmath}
\usepackage{amsfonts}
\usepackage{amssymb}
\usepackage{graphicx}
\usepackage{color}%
\providecommand{\U}[1]{\protect\rule{.1in}{.1in}}
\newtheorem{theorem}{Theorem}[section]

\newtheorem{lemma}[theorem]{Lemma}

\newtheorem{problem}{Problem}
\newtheorem{proposition}[theorem]{Proposition}

\newtheorem{remark}[theorem]{Remark}

\begin{document}

\title{\textbf{Broadcasts in Graphs:\ Diametrical Trees}\thanks{To appear in the
Australasian Journal of Combinatorics}}
\author{L. Gemmrich\thanks{Undergraduate research student}$\ ^{\ddag}$\ and C.M.
Mynhardt\thanks{Supported by the Natural Sciences and Engineering Research
Council of Canada.}\\Department of Mathematics and Statistics\\University of Victoria, Victoria, BC, \textsc{Canada}\\{\small lgemmrich@gmail.com}; {\small kieka@uvic.ca}}
\maketitle

\begin{abstract}
A dominating broadcast on a graph $G=(V,E)$ is a function $f:V\rightarrow
\{0,1,\dots,\operatorname{diam}(G)\}$ such that $f(v)\leq e(v)$ (the
eccentricity of $v$) for all $v\in V$, and each $u\in V$ is at distance at
most $f(v)$ from a vertex $v$ with $f(v)\geq1$. The upper broadcast domination
number of $G$ is $\Gamma_{b}(G)=\max\{\sum_{v\in V}f(v):f$ is a minimal
dominating broadcast on $G\}$. As shown by Erwin in [D.~Erwin, Cost domination
in graphs, Doctoral dissertation, Western Michigan University, 2001],
$\Gamma_{b}(G)\geq\operatorname{diam}(G)$ for any graph $G$.

We investigate trees whose upper broadcast domination number equal their
diameter and, among more general results, characterize caterpillars with this property.

\end{abstract}

\noindent\textbf{Keywords:\hspace{0.1in}}broadcast on a graph, dominating
broadcast, minimal dominating broadcast, upper broadcast domination number

\noindent\textbf{AMS 2010 Subject Classification Number:\hspace{0.1in}}05C69,
05C05, 05C12

\section{Introduction}

Suppose a telecommunications company has to provide radio coverage to a
collection of geographic regions. A single tower transmitting with a strength
(or cost) of one unit can provide coverage to the region it is located in and
all regions immediately adjacent to it. The company aims to minimize its
expenses by erecting as few towers as possible. If we consider each region as
a vertex of a graph $G$, where two vertices are adjacent if their
corresponding geographic regions are adjacent, then any \emph{dominating set}
$S$ (i.e. each vertex of $G$ belongs to $S$ or is adjacent to a vertex in $S$)
represents a suitable arrangement of radio towers, and a dominating set of
minimum cardinality represents a minimum cost arrangement. However, if the
company is able to build its towers with varying signal strength so that a
tower may transmit its signal a greater distance, but at a proportionally
greater cost, the total cost could be significantly less than for the former
arrangement. This situation can be modelled with a \textbf{broadcast} on $G$,
as defined below.

Unless stated otherwise, all graphs considered here are assumed to be simple,
nontrivial and connected. For undefined graph theoretic concepts and
terminology we refer the reader to \cite{CLZ} and \cite{HHS}.

A \emph{caterpillar} is a tree of order at least three, the removal of whose
leaves produces a path. We use standard notation for functions and write
$f:A\rightarrow B$ to denote the fact that $f$ is a function from $A$ to $B$;
we also write $f=\{(a,f(a)):a\in A\}$. If $f$ and $g$ are functions with the
same domain $A$ such that $g(a)\leq f(a)$ for each $a\in A$, we write $g\leq
f$. If in addition $g(a)<f(a)$ for at least one $a\in A$, we write $g<f$.

As usual we denote the domination and upper domination numbers of a graph $G$
by $\gamma(G)$ and $\Gamma(G)$, respectively. A \emph{broadcast} on a graph
$G=(V,E)$ is a function $f:V\rightarrow\{0,1,\dots,\operatorname{diam}(G)\}$
such that $f(v)\leq e(v)$ (the eccentricity of $v$) for all $v\in V$. A
broadcast $f$ on $G$ is \emph{dominating} if each $u\in V$ is at distance at
most $f(v)$ from a vertex $v$ with $f(v)\geq1$, and \emph{minimal dominating}
if no broadcast $f^{\prime}$ on $G$ with $f^{\prime}<f$ is dominating. The
\emph{cost} of a broadcast $f$ is $\sigma(f)=\sum_{v\in V}f(v)$. The
\emph{broadcast domination number} of $G$ is%
\[
\gamma_{b}(G)=\min\{\sigma(f):f\text{ is a dominating broadcast on }G\},
\]
and the \emph{upper broadcast domination number} of $G$ is%
\[
\Gamma_{b}(G)=\max\{\sigma(f):f\text{ is a minimal dominating broadcast on
}G\}.
\]
Broadcast domination was introduced by Erwin \cite{Ethesis, Epaper}, who
proved the bounds%
\begin{equation}
\gamma_{b}(G)\leq\min\{\gamma(G),\operatorname{rad}(G)\}\leq\max
\{\Gamma(G),\operatorname{diam}(G)\}\leq\Gamma_{b}(G) \label{eq_bound}%
\end{equation}
for any graph $G$. Graphs for which $\gamma_{b}(G)=\operatorname{rad}(G)$ are
called \emph{radial graphs}. Radial trees are characterized in \cite{Herke,
HM}. The upper broadcast domination number\emph{ }$\Gamma_{b}(G)$ is also
studied in \cite{Ahmadi, BF, Dunbar, MR}. Other studies of broadcast
domination can be found in \cite{BouchS, BZ, BMT, BS, CHM, DDH, HL, HS,
RadKhos, Scott, SM, MT, MW, MW2, SS}.

Our purpose is to investigate trees whose upper broadcast domination number
equals their diameter. Following the terminology for broadcast domination
numbers, we call such trees \emph{diametrical trees}. The characterization of
diametrical trees is listed as an open problem in \cite{MR}.

After presenting further definitions and known results in Section
\ref{Sec_Def}, we state a number of lemmas concerning properties of
non-diametrical trees in Section \ref{Sec_Non-D}. To avoid interrupting the
flow of the proof of our main theorem, we defer the proofs of all lemmas to
Section \ref{Sec_Proofs}. A consequence of these lemmas is that a tree
containing a path of length at least three, internally disjoint from a
diametrical path, is non-diametrical. This result hints that the caterpillars
may contain classes of diametrical trees, which is indeed the case. Our goal
is to prove the characterization of diametrical caterpillars stated in Theorem
\ref{Thm_Cater} below, which we do in Section~\ref{Sec_Cater}. We conclude
with open problems in Section \ref{Sec_Problems}.

\begin{theorem}
\label{Thm_Cater}A caterpillar $T$ with diametrical path $P:v_{0},v_{1}%
,\ldots,v_{d}$ is diametrical if and only if \vspace{-0.08in}

\begin{enumerate}
\item[$(i)$] each $v_{i},\ i\in\{1,\ldots,d-1\}$, is adjacent to at most two
leaves,\vspace{-0.08in}

\item[$(ii)$] for any $i\in\{1,\ldots,d-2\}$, $\min\{\deg_{T}(v_{i}),\deg
_{T}(v_{i+1})\}=2$,\vspace{-0.08in}

\item[$(iii)$] whenever $v_{i}$ and $v_{j},\ i<j$, are adjacent to at least
two leaves each, there exists an index $k,\ i<k<j$, such that $\deg_{T}%
(v_{k})=\deg_{T}(v_{k+1})=2$.
\end{enumerate}
\end{theorem}

\section{Definitions and Known Results\label{Sec_Def}}

For a broadcast $f$ on a graph $G=(V,E)$, define $V_{f}^{+}=\{v\in
V:f(v)>0\}$. The vertices in $V_{f}^{+}$ are called \emph{broadcast vertices}.
A vertex $u$ \emph{hears} the broadcast $f$ from some vertex $v\in V_{f}^{+}$,
and $v$ $f$-\emph{dominates} $u$, if the distance $d(u,v)\leq f(v)$. An edge
$uw$ \emph{hears} $f$ if both $u$ and $w$ hear $f$ from the same vertex $v\in
V_{f}^{+}$. A vertex $v\in V_{f}^{+}$ \emph{overdominates }a vertex $u$ if
$d(u,v)<f(v)$. For $v\in V_{f}^{+}$, define the

\begin{itemize}
\item $f$-\emph{neighbourhood} of $v$ as $N_{f}[v]=\{u\in V(G):d(u,v)\leq
f(v)\},$

\item $f$-\emph{boundary} of $v$ as\emph{ }$B_{f}(v)=\{u\in
V(G):d(u,v)=f(v)\},$

\item $f$-\emph{private neighbourhood }of\emph{ }$v$ as $\operatorname{PN}%
_{f}(v)=\{u\in N_{f}[v]:u\notin N_{f}[w]$ for all $w\in V^{+}-\{v\}\}$,

\item $f$-\emph{private boundary} of $v$ as $\operatorname{PB}_{f}(v)=\{u\in
N_{f}[v]:u$ is not dominated by $(f-\{(v,f(v))\})\cup\{(v,f(v)-1)\}\}$.
\end{itemize}

\begin{figure}[pb]%
\centering
\includegraphics[
height=1.0533in,
width=3.0441in
]%
{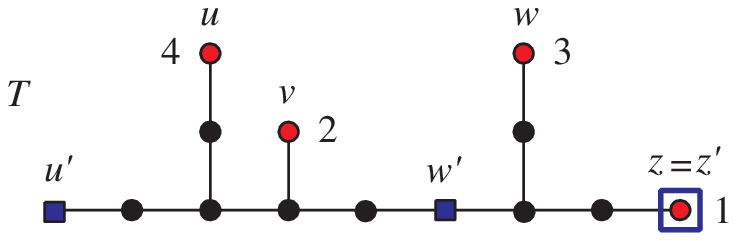}%
\caption{A tree $T$ with a dominating broadcast $f$ such that
$\operatorname{PB}_{f}(x)=\{x^{\prime}\}$ for each $x\in\{u,w,z\}$, and
$\operatorname{PB}_{f}(v)=\varnothing.$}%
\label{Fig_Ex1}%
\end{figure}

Note that if $f(v)=1$, then $\operatorname{PB}_{f}(v)=\operatorname{PN}%
_{f}(v)$, and if $f(v)\geq2$, then $\operatorname{PB}_{f}(v)=B_{f}%
(v)\cap\operatorname{PN}_{f}(v)$. For example, consider the tree $T$ in Figure
\ref{Fig_Ex1}. The broadcast $f$ defined by
$f(u)=4,\ f(v)=2,\ f(w)=3,\ f(z)=1$ and $f(x)=0$ otherwise is a dominating
broadcast such that $\operatorname{PB}_{f}(x)=\{x^{\prime}\}$ for each
$x\in\{u,w,z\}$, and $\operatorname{PB}_{f}(v)=\varnothing$.%

The property that makes a dominating broadcast minimal dominating, determined
in \cite{Ethesis} and stated in \cite{MR} in terms of private
boundaries\emph{,} is essential in the study of upper broadcast numbers. We
state it again here.

\begin{proposition}
\label{PropMinimal}\emph{\cite{Ethesis}}\hspace{0.1in}A dominating broadcast
$f$ is a minimal dominating broadcast if and only if $\operatorname{PB}%
_{f}(v)\neq\varnothing$ for each $v\in V_{f}^{+}$.
\end{proposition}

\begin{figure}[ptb]%
\centering
\includegraphics[
height=1.0533in,
width=3.0441in
]%
{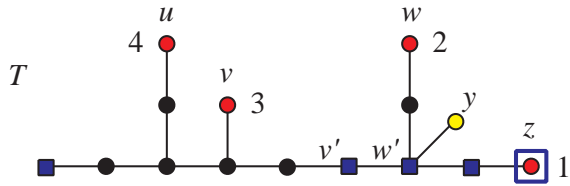}%
\caption{A tree $T$ with a broadcast $f$ such that $\operatorname{PB}%
_{f}(x)\neq\varnothing$ for each $x\in V_{f}^{+}$, but $f$ cannot be extended
to a minimal dominating broadcast.}%
\label{Fig_Ex2}%
\end{figure}

By Proposition \ref{PropMinimal} the broadcast $f$ in Figure \ref{Fig_Ex1},
although dominating, is not minimal dominating. The broadcast $f^{\prime
}=(f-\{(v,2)\})\cup\{(v,0)\}$ is a minimal dominating broadcast on $T$. In
general it is not true that if $f$ is a dominating broadcast on a graph $G$,
then some broadcast $f^{\prime}$ with $f^{\prime}\leq f$ is a minimal
dominating broadcast on $G$, nor is it necessarily true that if $f$ is a
broadcast on $G$ such that $\operatorname{PB}_{f}(v)\neq\varnothing$ for each
$v\in V_{f}^{+}$, then some broadcast $f^{\prime}$ with $f\leq f^{\prime}$ is
a minimal dominating broadcast on $G$. Consider the tree $T$ and broadcast $f$
shown in Figure \ref{Fig_Ex2}. Here, $\operatorname{PB}_{f}(x)=\{x^{\prime}\}$
for each $x\in\{v,w\}$ and $y$ is not $f$-dominated. Moreover, $f$ cannot be
extended to a broadcast that dominates $y$ without leaving $v$ or $w$ with an
empty private boundary.%

It is well known that any independent set of vertices in a graph $G$ can be
extended to a maximal (but not necessarily maximum) independent set of $G$,
and that a maximal independent set is also a minimal dominating set (cf.
\cite[pp.~70 -- 71]{HHS}). Denoting the cardinality of a maximum independent
set of $G$ by $\alpha(G)$, it follows that $\alpha(G)\leq\Gamma(G)$ for all
graphs $G$.

\begin{remark}
\label{Rem_Bound}\emph{\cite{Ethesis}}\hspace{0.1in}The characteristic
function of a minimal dominating set in a graph $G$ is a minimal dominating
broadcast on $G$. Hence $\Gamma_{b}(G)\geq\Gamma(G)\geq\alpha(G)$ for any
graph $G$.
\end{remark}

\begin{proposition}
\label{Prop_Disjoint}\emph{\cite{Ethesis}}\hspace{0.1in}If $f$ is a broadcast
on a graph $G$ and for each $i\in\{1,2\}$ we have $u_{i}\in V_{f}^{+}$,
$u_{i}^{\prime}\in\operatorname{PB}_{f}(u_{i})$,\ where $u_{1}\neq u_{2}$, and
$P_{i}$ is a $u_{i}$ -- $u_{i}^{\prime}$ geodesic, then $P_{1}$ and $P_{2}$
are disjoint.
\end{proposition}

Using Proposition \ref{Prop_Disjoint}, Erwin \cite{Ethesis} shows that
$\Gamma_{b}(G)\leq|E(G)|$ for any graph $G$, and together with the lower bound
(\ref{eq_bound}) this implies that $\Gamma_{b}(P_{n})=n-1$ for each $n\geq2$.
Proposition \ref{Prop_Disjoint} is used frequently in the proofs in Section
\ref{Sec_Proofs}.

\section{Non-Diametrical Trees\label{Sec_Non-D}}

\begin{figure}[ptb]%
\centering
\includegraphics[
height=1.0101in,
width=2.8599in
]%
{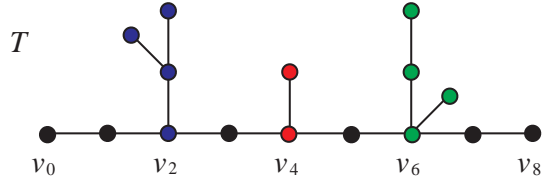}%
\caption{A tree $T$ and diametrical path $P:v_{0},v_{1},...,v_{8}$ such that
$T_{2}\protect\cong K_{1,3},\ T_{4}\protect\cong K_{2},\ T_{6}\protect\cong
P_{4}$ and $T_{i}=K_{1}$ otherwise.}%
\label{Fig_Ex3}%
\end{figure}

In this section we state a number of sufficient conditions for a tree $T$ to
be non-diametrical. The proofs are given in Section \ref{Sec_Proofs}. We
assume throughout that $T$ has diameter $d$ and a diametrical path
$P:v_{0},v_{1},\ldots,v_{d}$. For each $i\in\{0,\ldots,d\}$, let $T_{i}$ be
the subtree of $T$ induced by all vertices that are connected to $v_{i}$ by
paths that are internally disjoint from $P$. Note that $T_{i}=K_{1}$ if and
only if $i\in\{0,d\}$, or $i\in\{1,\ldots,d-1\}$ and $\deg(v_{i})=2$. For
example, in the tree $T$ in Figure \ref{Fig_Ex3}, $T_{2}\cong K_{1,3}%
,\ T_{4}\cong K_{2},\ T_{6}\cong P_{4}$ and $T_{i}=K_{1}$ for $i\in
\{0,1,3,5,7,8\}$.

A \emph{stem} of a tree $T\ncong K_{2}$ is a vertex adjacent to a leaf and a
\emph{strong stem} is a stem that is adjacent to at least two leaves; in
Figure \ref{Fig_Ex3}, $v_{1},v_{4}$ and $v_{6}$ are (not the only) examples of
stems. The complete bipartite graph $K_{1,t},\ t\geq1$, is also called a
\emph{star. }Thus a tree\emph{ }$T$ with diametrical path $P$ as above is a
caterpillar\emph{ }if each $T_{i}$ is either a star or $K_{1}$.

\begin{lemma}
\label{Lem1-1}Let $T$ be a tree with diameter $d\geq3$ and diametrical path
$P:v_{0},v_{1},\ldots,v_{d}$. If there exists an $i\in\{1,\ldots,d-2\}$ such
that each of $v_{i}$ and $v_{i+1}$ is adjacent to a leaf other than $v_{0}$
(if $i=1$) or $v_{d}$ (if $i+1=d-1$), then $\Gamma_{b}(T)>\operatorname{diam}%
(T)$.
\end{lemma}

%

\begin{lemma}
\label{Lem_alpha_Ti}If there exists a subscript $i\in\{2,\ldots,d-2\}$ such
that $T_{i}$ has an independent set of cardinality $3$ that dominates but does
not contain $v_{i}$, or if $\max\{\deg(v_{1}),\deg(v_{d-1})\}\geq4$, then
$\Gamma_{b}(T)>\operatorname{diam}(T)$.
\end{lemma}

\begin{lemma}
\label{Lem_alpha_2}If there exists a subscript $i\in\{2,\ldots,d-2\}$ such
that $T_{i}$ has an independent set of cardinality $2$ that does not dominate
$v_{i}$, then $\Gamma_{b}(T)>\operatorname{diam}(T)$.
\end{lemma}

\begin{lemma}
\label{Lem_diam3}If $\operatorname{diam}(T_{i})\geq4$ for some $i$, or if
$\operatorname{diam}(T_{i})=3$ and $v_{i}$ is a peripheral vertex of $T_{i}$,
then $\Gamma_{b}(T)>\operatorname{diam}(T)$.
\end{lemma}

By Lemmas \ref{Lem_alpha_Ti} -- \ref{Lem_diam3}, if $T$ is a diametrical tree,
then each $T_{i}$ is isomorphic to either $K_{1}$, $K_{2}$, $P_{3}$ with
$v_{i}$ either a leaf or the stem of $P_{3}$, or $P_{4}$ with $v_{i}$ being a
stem of $P_{4}$. Thus, diametrical trees are \textquotedblleft
nearly\textquotedblright\ caterpillars. We henceforth restrict our
investigation to caterpillars. By Lemma \ref{Lem1-1}, if $T_{i}\cong K_{2}$,
we may assume that neither $T_{i-1}$ nor $T_{i+1}$ is isomorphic to $K_{2}$.
If $T_{i}\cong P_{3}$ with $v_{i}$ being a leaf of $P_{3}$, or if $T_{i}\cong
P_{4}$, then $T$ is not a caterpillar and we ignore these cases. We give one
more sufficient condition for a caterpillar to be non-diametrical.

\begin{lemma}
\label{Lem_K1,2_Caterpillar}Let $T$ be a caterpillar with diametrical path
$P:v_{0},v_{1},\ldots,v_{d}$. If two\ vertices $v_{i},v_{i+2k}$ are strong
stems, for some $i\geq1$ and some integer $k$ such that $i+2k\leq d-1$, and
$v_{i+2r}$ is a stem for each $r\in\{1,\ldots,k-1\}$, then $\Gamma_{b}(T)>d$.
\end{lemma}

\section{Diametrical Caterpillars\label{Sec_Cater}}

If $T$ is a diametrical caterpillar, then $T$ does not satisfy the hypothesis
of any of Lemmas \ref{Lem1-1} -- \ref{Lem_K1,2_Caterpillar}. In this section
we show that the converse is also true:\ If the caterpillar $T$ does not
satisfy the hypothesis of any of Lemmas \ref{Lem1-1} --
\ref{Lem_K1,2_Caterpillar}, then $T$ is diametrical. The negation of these
hypotheses, applied to caterpillars, gives the characterization of diametrical
caterpillars stated in Theorem \ref{Thm_Cater}, which we restate here for convenience.

\bigskip

\noindent\textbf{Theorem \ref{Thm_Cater}\hspace{0.1in}}\emph{A caterpillar
}$T$\emph{ with diametrical path }$P:v_{0},v_{1},\ldots,v_{d}$\emph{ is
diametrical if and only if\vspace{-0.08in}}

\begin{enumerate}
\item[$(i)$] \emph{each }$v_{i},\ i\in\{1,\ldots,d-1\}$\emph{, is adjacent to
at most two leaves},\emph{\vspace{-0.08in}}

\item[$(ii)$] \emph{for any }$i\in\{1,\ldots,d-2\}$\emph{, }$\min\{\deg
_{T}(v_{i}),\deg_{T}(v_{i+1})\}=2$\emph{,\vspace{-0.08in}}

\item[$(iii)$] \emph{whenever }$v_{i}$\emph{ and }$v_{j},\ i<j$\emph{, are
strong stems, there exists an index }$k,\ i<k<j$\emph{, such that }$\deg
_{T}(v_{k})=\deg_{T}(v_{k+1})=2$\emph{.}
\end{enumerate}

\noindent\textbf{Proof.\hspace{0.1in}}\label{here2}Suppose $T$ is a
diametrical caterpillar. By Lemma \ref{Lem_alpha_Ti}, each $v_{i}%
,\ i\in\{2,\ldots,d-2\}$, is adjacent to at most two leaves, while $v_{1}$ and
$v_{d-1}$ are adjacent to at most one leaf other than $v_{0}$ and $v_{d}$,
respectively, hence $(i)$ holds. Similarly, condition $(ii)$ follows directly
from Lemma \ref{Lem1-1}. For $(iii)$, condition $(ii)$ already implies that of
any two consecutive internal vertices of $P$, at least one has degree $2$.
Lemma \ref{Lem_K1,2_Caterpillar} now implies that if $v_{i}$ and $v_{j}$ are
both strong stems, then some pair of consecutive strong stems between $v_{i}$
and $v_{j}$ (inclusive) are separated by at least two vertices of degree $2$.
Hence $(iii)$ holds.

For the converse, note that the only caterpillars of diameter three or less
that satisfy conditions $(i)$ -- $(iii)$ are $P_{3},\ P_{4}$ and the tree
obtained by joining a new leaf to a stem of $P_{4}$. It is easy to verify that
they are diametric. Assume that Theorem \ref{Thm_Cater} is false and let $T$
be a smallest non-diametrical caterpillar that satisfies $(i)$ -- $(iii)$.
Then $T$ has diameter at least four. We state two more lemmas, the proofs of
which are also given in Section \ref{Sec_Proofs}.

\begin{lemma}
\label{Lem_non-diam3}No vertex of $T$ is a strong stem.
\end{lemma}

\begin{lemma}
\label{Lem_K2}No vertex $v_{i},\ i\in\{2,\ldots,d-2\}$, is adjacent to a leaf.
\end{lemma}

By Lemmas \ref{Lem_non-diam3} and \ref{Lem_K2}, $T=P_{d+1}$, which is
impossible because Erwin \cite{Ethesis} showed that $\Gamma_{b}(P_{n}%
)=n-1=\operatorname{diam}(P_{n})$ for all $n\geq2$.~$\blacksquare$

\section{Proofs of Lemmas\label{Sec_Proofs}}

This section contains the proofs of Lemmas \ref{Lem1-1} -- \ref{Lem_K2},
restated here for convenience.\medskip

\noindent\textbf{Lemma \ref{Lem1-1}\hspace{0.1in}}\emph{Let }$T$\emph{ be a
tree with diameter }$d\geq3$\emph{ and diametrical path }$P:v_{0},v_{1}%
,\ldots,v_{d}$\emph{. If some }$v_{i}$\emph{ and }$v_{i+1}$\emph{, }%
$i\in\{1,\ldots,d-2\}$\emph{, are adjacent to leaves other than }$v_{0}$\emph{
or }$v_{d}$\emph{, then }$\Gamma_{b}(T)>\operatorname{diam}(T)$\emph{.}

\medskip

\noindent\textbf{Proof.\hspace{0.1in}}Suppose the hypothesis of the lemma is
satisfied. Say $v_{i}$ is adjacent to the leaf $\ell$ and $v_{i+1}$ is
adjacent to the leaf $\ell^{\prime}$. Define the broadcast $g$ by
$g(v_{0})=i+1$, $g(v_{d})=d-i$ and $g(x)=0$ otherwise. Then $\ell
\in\operatorname{PB}_{g}(v_{0})$ and $\ell^{\prime}\in\operatorname{PB}%
_{g}(v_{d})$, hence $\operatorname{PB}_{g}(x)\neq\varnothing$ for all $x\in
V_{g}^{+}$. If $g$ is also dominating, let $f=g$; otherwise, let $T^{\prime}$
be the subgraph of $T$ induced by all vertices that are not $g$-dominated, let
$S$ be a maximal independent set of $T^{\prime}$ and define the broadcast $f$
by $f(x)=g(x)$ if $x\in V(T)-V(T^{\prime})$, $f(x)=1$ if $x\in S$ and $f(x)=0$
if $x\in V(T^{\prime})-S$. By definition, $f$ is a dominating broadcast on
$T$. Since $\ell$ and $\ell^{\prime}$ are leaves, no vertex in $S$ is adjacent
to $\ell$ or $\ell^{\prime}$, hence $\ell\in\operatorname{PB}_{f}(v_{0})$ and
$\ell^{\prime}\in\operatorname{PB}_{f}(v_{d})$. Since no vertex in $S$ hears
the broadcast $g$, $x\in\operatorname{PB}_{f}(x)$ for each $x\in S$. Hence, by
Proposition \ref{PropMinimal}, $f$ is a minimal dominating broadcast.
Moreover, $\sigma(f)\geq i+1+d-i=d+1$ and the result follows.~$\blacksquare
$\medskip

The proof of the next lemma is illustrated in Figure \ref{Fig_alpha_Ti}%
.\medskip

\noindent\textbf{Lemma \ref{Lem_alpha_Ti}\hspace{0.1in}}\emph{If there exists
a subscript }$i\in\{2,\ldots,d-2\}$\emph{ such that }$T_{i}$\emph{ has an
independent set of cardinality }$3$\emph{ that dominates }$T_{i}$\emph{ but
does not contain }$v_{i}$\emph{, or if }$\max\{\alpha(T_{1}),\alpha
(T_{d-1})\}\geq2$\emph{, then }$\Gamma_{b}(T)>\operatorname{diam}(T)$%
\emph{.}\medskip

\noindent\textbf{Proof.\hspace{0.1in}}We may assume that $T$ does not satisfy
the hypothesis of Lemma \ref{Lem1-1}, otherwise we are done.\textbf{ }Suppose
$\alpha(T_{1})=t\geq2$. See Figure \ref{Fig_alpha_Ti}(a). Since $v_{0}$ is a
peripheral vertex of $T$, no vertex of $T_{1}$ is at distance greater than one
from $v_{1}$. Hence $T_{1}=K_{1,t}$ and, by Lemma \ref{Lem1-1}, $v_{2}$ is not
adjacent to a leaf. Let $S$ be the set consisting of $v_{0}$ and the $t$
leaves of $T_{1}$, and define the broadcast $g$ by $g(v_{d})=d-2$, $d(x)=1$ if
$x\in S$ and $g(x)=0$ otherwise. Then $v_{2}\in\operatorname{PB}_{g}(v_{d})$
and $x\in\operatorname{PB}_{g}(x)$ for each $x\in S$, hence $\operatorname{PB}%
_{g}(x)\neq\varnothing$ for all $x\in V_{g}^{+}$. If $g$ is dominating, let
$f=g$, otherwise let $T^{\prime}$ be the subgraph of $T$ induced by all
vertices that are not dominated by $g$. Since $v_{2}$ is not adjacent to a
leaf, there exists a maximal independent set $X$ of $T^{\prime}$ that does not
contain a vertex adjacent to $v_{2}$. Define the broadcast $f$ by $f(x)=g(x)$
if $x\in V(T)-V(T^{\prime})$, $f(x)=1$ if $x\in X$ and $f(x)=0$ if $x\in
V(T^{\prime})-X$. Then $v_{2}\in\operatorname{PB}_{f}(v_{d})$ and
$x\in\operatorname{PB}_{f}(x)$ for each $x\in S\cup X$, so $f$ is a minimal
dominating broadcast on $T$ with $\sigma(f)\geq t+1+d-2>d$. Hence $\Gamma
_{b}(T)>d$.

If $\alpha(T_{d-1})\geq2$ the result follows similarly. Hence assume some
$T_{i},\ i\in\{2,\ldots,d-2\}$, has an independent set of cardinality $3$ that
dominates $T_{i}$ but does not contain $v_{i}$. Then $T_{i}$ has a maximal
independent set $S$ of cardinality $c\geq3$ such that $v_{i}\notin S$. Define
the broadcast $g$ by $g(v_{0})=i-1$, $g(v_{d})=d-i-1$, $g(x)=1$ if $x\in S$
and $g(x)=0$ otherwise. Since $v_{i}\notin S$, $v_{i-1}\in\operatorname{PB}%
_{g}(v_{0})$ and $v_{i+1}\in\operatorname{PB}_{g}(v_{d})$. In addition,
$x\in\operatorname{PB}_{g}(x)$ for each $x\in S$. If $i\geq3$ and $v_{i-1}$ is
adjacent to a leaf, then we may assume, by Lemma \ref{Lem1-1}, that $v_{i-2}$
is not adjacent to a leaf (other than $v_{0}$ if $i=3$). Similarly, if $i\leq
d-3$ and $v_{i+1}$ is adjacent to a leaf, we may assume that $v_{i+2}$ is not
adjacent to a leaf (other than $v_{d}$ if $i=d-3)$. Let $T^{\prime}$ be the
subgraph of $T$ induced by the vertices that are not dominated by $g$ and
choose a maximal independent set $X$ of $T^{\prime}$ as follows.

\begin{itemize}
\item If $T^{\prime}$ has a maximal independent set that does not contain a
vertex adjacent to $v_{i-1}$ or to $v_{i+1}$, let $X$ be such a set. See
Figure \ref{Fig_alpha_Ti}(b).%

\begin{figure}[ptb]%
\centering
\includegraphics[
height=3.0364in,
width=4.4927in
]%
{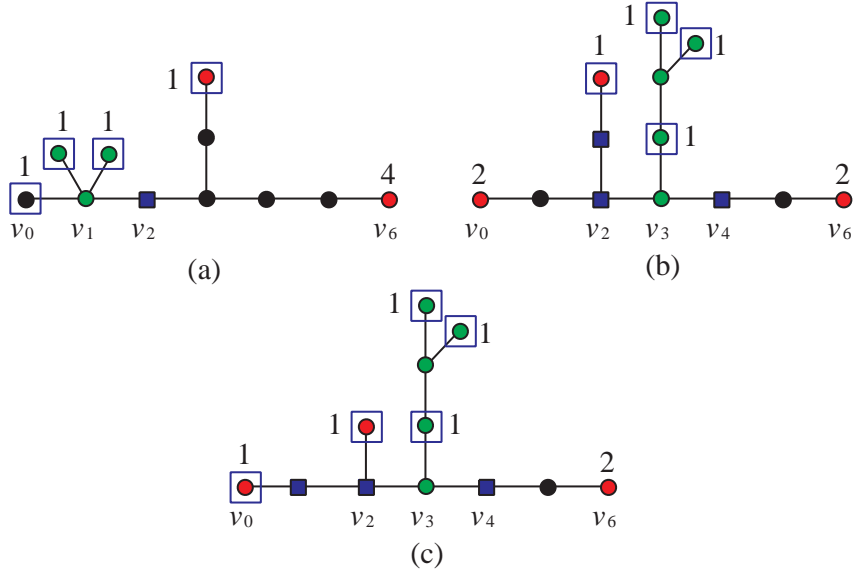}%
\caption{An illustration of the proof of Lemma \ref{Lem_alpha_Ti}.}%
\label{Fig_alpha_Ti}%
\end{figure}

\item If each maximal independent set of $T^{\prime}$ contains a vertex
adjacent to $v_{i-1}$ (or $v_{i+1}$ or both), then $v_{i-1}$ (or $v_{i+1}$) is
adjacent to a leaf. Then $v_{i-2}$ (or $v_{i+2}$) is not adjacent to a leaf,
and there exists a maximal independent set of $T^{\prime}$ that contains no
vertex adjacent to $v_{i-2}$ (or $v_{i+2}$); let $X$ be such a set. See Figure
\ref{Fig_alpha_Ti}(c).
\end{itemize}

Define the broadcast $f$ on $T$ as follows. If neither $v_{i-1}$ nor $v_{i+1}$
is adjacent to a leaf, let%
\[
f(x)=\left\{
\begin{tabular}
[c]{ll}%
$g(x)$ & if $x\in V(T)-V(T^{\prime})$\\
$1$ & if $x\in X$\\
$0$ & otherwise.
\end{tabular}
\right.
\]
Then $f$ is a dominating broadcast such that $\sigma(f)\geq i-1+d-i-1+c>d$,
$v_{i-1}\in\operatorname{PB}_{f}(v_{0})$, $v_{i+1}\in\operatorname{PB}%
_{f}(v_{d})$ and $x\in\operatorname{PB}_{f}(x)$ for each $x\in S\cup X$.

If $v_{i-1}$ is adjacent to a leaf and $v_{i+1}$ is not, let%
\[
f(x)=\left\{
\begin{tabular}
[c]{ll}%
$i-2$ & if $x=v_{0}$\\
$g(x)$ & if $x\in V(T-v_{0})-V(T^{\prime})$\\
$1$ & if $x\in X$\\
$0$ & otherwise.
\end{tabular}
\ \ \right.
\]
Then $|X|\geq1$ and $v_{i-1}$ hears $f$ from an adjacent leaf. Hence $f$ is a
dominating broadcast such that $\sigma(f)\geq i-2+d-i-1+c+|X|>d,\ v_{i-2}%
\in\operatorname{PB}_{f}(v_{0}),\ v_{i+1}\in\operatorname{PB}_{f}(v_{d})$ and
$x\in\operatorname{PB}_{f}(x)$ for each $x\in S\cup X$.

Similarly, if $v_{i+1}$ is adjacent to a leaf and $v_{i-1}$ is not, let%
\[
f(x)=\left\{
\begin{tabular}
[c]{ll}%
$d-i-2$ & if $x=v_{d}$\\
$g(x)$ & if $x\in V(T-v_{d})-V(T^{\prime})$\\
$1$ & if $x\in X$\\
$0$ & otherwise.
\end{tabular}
\right.
\]

Finally, if both $v_{i-1}$ and $v_{i+1}$ are adjacent to leaves, define $f$ by%
\[
f(x)=\left\{
\begin{tabular}
[c]{ll}%
$i-2$ & if $x=v_{0}$\\
$d-i-2$ & if $x=v_{d}$\\
$g(x)$ & if $x\in V(T-\{v_{0},v_{d}\})-V(T^{\prime})$\\
$1$ & if $x\in X$\\
$0$ & otherwise.
\end{tabular}
\ \right.
\]
Now $|X|\geq2$ and $f$ is a dominating broadcast such that $\sigma(f)\geq
i-2+d-i-2+c+|X|\geq d-4+c+|X|>d$, $v_{i-2}\in\operatorname{PB}_{f}%
(v_{0}),\ v_{i+2}\in\operatorname{PB}_{f}(v_{d})$ and $x\in\operatorname{PB}%
_{f}(x)$ if $x\in S\cup X$.

Hence in each case $f$ is a minimal dominating broadcast such that
$\sigma(f)>d$, which implies that $\Gamma_{b}(T)>\operatorname{diam}%
(T).~\blacksquare\medskip$

\noindent\textbf{Lemma \ref{Lem_alpha_2}\hspace{0.1in}}\emph{If there exists a
subscript }$i\in\{2,\ldots,d-2\}$\emph{ such that }$T_{i}$\emph{ has an
independent set of cardinality }$2$\emph{ that does not dominate }$v_{i}%
$\emph{, then }$\Gamma_{b}(T)>\operatorname{diam}(T)$\emph{.}$\medskip$

\noindent\textbf{Proof.\hspace{0.1in}}Suppose $T_{i}$ has an independent set
$D$ of cardinality $2$ that does not dominate $v_{i}$. If every maximal
independent set of $T_{i}$ that contains $D$, but not $v_{i}$, dominates
$v_{i}$, the result follows from Lemma \ref{Lem_alpha_Ti}. Hence assume this
is not the case (in particular, $v_{i}$ is not a stem) and let $S$ be a
maximal independent set of cardinality $c\geq2$ of $T_{i}-v_{i}$ containing no
vertex adjacent to $v_{i}$. Define the broadcast $g$ on $T$ by $g(v_{0}%
)=i,\ g(v_{d})=d-i-1,\ g(x)=1$ for each $x\in S$ and $g(x)=0$ otherwise. Note
that $v_{i}\in\operatorname{PB}_{g}(v_{0})$, $v_{i+1}\in\operatorname{PB}%
_{g}(v_{d}),\ x\in\operatorname{PB}_{g}(x)$ for each $x\in X$ and
$\sigma(g)\geq i+d-i-1+c>d$. We can now proceed as in the proof of Lemma
\ref{Lem_alpha_Ti} to construct a minimal dominating broadcast $f$ on $T$ such
that $\sigma(f)\geq\sigma(g)>d$ to obtain that $\Gamma_{b}(T)>d$. The details
are omitted.$~\blacksquare\medskip$

\noindent\textbf{Lemma \ref{Lem_diam3}\hspace{0.1in}}\emph{If }%
$\operatorname{diam}(T_{i})\geq4$\emph{ for some }$i$\emph{, or if
}$\operatorname{diam}(T_{i})=3$\emph{ and }$v_{i}$\emph{ is a peripheral
vertex of }$T_{i}$\emph{, then }$\Gamma_{b}(T)>\operatorname{diam}(T)$%
\emph{.}$\medskip$

\noindent\textbf{Proof.\hspace{0.1in}}If $\operatorname{diam}(T_{i})\geq5$,
then $T_{i}$ contains a subgraph isomorphic to $P_{6}$, which, regardless of
which vertex of $P_{6}$ corresponds to $v_{i}$, has an independent set of
cardinality $3$ that dominates but does not contain $v_{i}$, and the result
follows from Lemma \ref{Lem_alpha_Ti}. If $\operatorname{diam}(T_{i})=4$ and
$v_{i}$ corresponds to a stem of a subgraph isomorphic to $P_{5}$, the result
follows similarly.

Suppose $\operatorname{diam}(T_{i})=k\in\{3,4\}$ and $v_{i}$ is a peripheral
vertex of $T_{i}$. Then $v_{i}$ is not a stem. Let $\ell$ be a vertex of
$T_{i}$ at distance $k$ from $v_{i}$. Define the broadcast $g$ on $T$ by
$g(\ell)=k,\ g(v_{0})=i-1,\ g(v_{d})=d-i-1$ and $g(x)=0$ otherwise. Then
$v_{i}\in\operatorname{PB}_{g}(\ell),$ $v_{i-1}\in\operatorname{PB}_{g}%
(v_{0})$ and $v_{i+1}\in\operatorname{PB}_{g}(v_{d})$, while $\sigma
(g)=i-1+d-i-1+k>d$. Possibly $v_{i-1}$ or $v_{i+1}$ is a stem, or both are. We
proceed as in the proof of Lemma \ref{Lem_alpha_Ti} to show that $\Gamma
_{b}(T)>d$.

Finally, suppose $\operatorname{diam}(T_{i})=4$ and $v_{i}$ is the central
vertex of a subgraph $H\cong P_{5}$ of $T_{i}$. Let $\ell_{1}$ and $\ell_{2}$
be the leaves of $H$ and let $w$ be the stem of $H$ adjacent to $\ell_{2}$. If
$v_{i}$ is a stem of $T_{i}$ the result again follows from Lemma
\ref{Lem_alpha_Ti}, hence assume $v_{i}$ is not a stem. Define the broadcast
$g$ by $g(\ell_{1})=2,\ g(\ell_{2})=1,\ g(v_{0})=i-1,\ g(v_{d})=d-i-1$ and
$g(x)=0$ otherwise. Then $v_{i}\in\operatorname{PB}_{g}(\ell_{1}%
),\ w\in\operatorname{PB}_{g}(\ell_{2}),$ $v_{i-1}\in\operatorname{PB}%
_{g}(v_{0})$ and $v_{i+1}\in\operatorname{PB}_{g}(v_{d})$, while
$\sigma(g)=i-1+d-i-1+3>d$. As before it (eventually) follows that $\Gamma
_{b}(T)>d$.~$\blacksquare\medskip$

\noindent\textbf{Lemma \ref{Lem_K1,2_Caterpillar}\hspace{0.1in}}\emph{Let }%
$T$\emph{ be a caterpillar with diametrical path }$P:v_{0},v_{1},\ldots,v_{d}%
$\emph{. If two\ vertices }$v_{i},v_{i+2k}$\emph{ are strong stems, for some
}$i\geq1$\emph{ and some integer }$k$\emph{ such that }$i+2k\leq d-1$\emph{,
and }$v_{i+2r}$\emph{ is a stem for each }$r\in\{1,\ldots,k-1\}$\emph{, then
}$\Gamma_{b}(T)>d$\emph{.}\medskip

\noindent\textbf{Proof.\hspace{0.1in}}Let $S$ be the set of leaves adjacent to
$v_{i+2t},\ t\in\{0,1,\ldots,k\}$, and $X=\{v_{i+1},v_{i+3},\ldots,$
$v_{i+2k-1}\}$. Then $S\cup X$ is independent. By the hypothesis, $|S|\geq
k+3$ and so $|S\cup X|\geq2k+3$. By Lemma \ref{Lem1-1} we may assume that
$\deg_{T}(x)=2$ for each $x\in X$, otherwise the result follows.

If $i=1$ and $i+2k=d-1$, then $S\cup X$ is a maximal independent set of $T$ of
cardinality at least $d+1$. Let $f$ be the characteristic function of $S\cup
X$.

If $i=1$ and $i+2k<d-1$, define the broadcast $f$ on $T$ by $f(x)=1$ if $x\in
S\cup X$, $f(v_{d})=d-i-2k-1$ and $f(x)=0$ otherwise. Then $x\in
\operatorname{PB}_{f}(x)$ for each $x\in S\cup X$ and $v_{i+2k+1}%
\in\operatorname{PB}_{f}(v_{d})$. Since $v_{i+2k}$ is a stem, we may assume
that $\deg(v_{i+2k+1})=2$, otherwise the result holds by Lemma \ref{Lem1-1}.
Therefore $f$ is a dominating broadcast, thus a minimal dominating broadcast,
and $\sigma(f)=|S\cup X|+d-i-2k-1\geq d+1$.

\begin{figure}[ptb]%
\centering
\includegraphics[
height=0.838in,
width=3.3615in
]%
{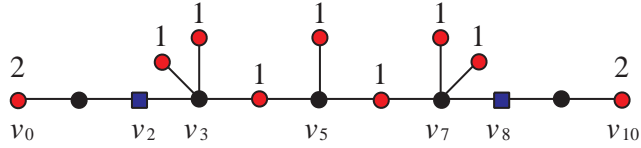}%
\caption{An illustration of the proof of Lemma \ref{Lem_K1,2_Caterpillar}.}%
\label{Fig_Cater}%
\end{figure}

If $i>1$ and $i+2k=d-1$, reverse the direction of $P$ and proceed as above.
Hence assume $1<i<i+2k<d-1$. See Figure \ref{Fig_Cater}, where $d=10,\ i=3$
and $k=2$. As above we may assume that $\deg(v_{i-1})=\deg(v_{i+2k+1})=2$.
Define the broadcast $f$ by $f(v_{0})=i-1,\ f(v_{d})=d-i-2k-1,\ f(x)=1$ for
each $x\in S\cup X$ and $f(x)=0$ otherwise. Then $f$ is a dominating broadcast
such that $\sigma(f)\geq d-2k-2+2k+3>d,\ v_{i-1}\in\operatorname{PB}_{f}%
(v_{0}),\ v_{i+2k+1}\in\operatorname{PB}_{f}(v_{d})$ and $x\in
\operatorname{PB}_{f}(x)$ for each $x\in S\cup X$. Hence $f$ is a minimal
dominating broadcast of $T$ such that $\sigma(f)>d$. The result now
follows.~$\blacksquare$\medskip%

Before proving Lemmas \ref{Lem_non-diam3} and \ref{Lem_K2} we state and prove
two additional lemmas. If $f$ is a broadcast on $T$ and $T^{\prime}$ is a
subtree of $T$, we define the \emph{restriction of }$f$\emph{ to} $T^{\prime}$
to be the broadcast $f^{\prime}=f\upharpoonleft T^{\prime}$ with
$V_{f^{\prime}}^{+}=V_{f}^{+}\cap V(T^{\prime})$ and $f^{\prime}(x)=f(x)$ for
all $x\in V(T^{\prime})$.

\begin{lemma}
\label{Lem_non-diam1}Suppose $T$ is a smallest non-diametrical caterpillar
that satisfies Theorem $\ref{Thm_Cater}(i)$ -- $(iii)$. Let $f$ be a minimal
dominating broadcast on $T$ such that $\sigma(f)>\operatorname{diam}(T)$. Then
$v_{0}\in V_{f}^{+}$ or $\{v_{0}\}=\operatorname{PB}_{f}(x)$ for some $x\in
V_{f}^{+}$, and a similar result holds for $v_{d}$.
\end{lemma}

\noindent\textbf{Proof.\hspace{0.1in}}Suppose the conclusion is false and say
$u\in V_{f}^{+}$ broadcasts to $v_{0}$, where $u\neq v_{0}$. Since
$\{v_{0}\}\neq\operatorname{PB}_{f}(u)$, there exists $b\in\operatorname{PB}%
_{f}(u)-\{v_{0}\}$. Possibly $b$ is a leaf adjacent to $v_{1}$, in which case
$v_{0}\in\operatorname{PB}_{f}(u),$ $\operatorname{diam}%
(T-b)=\operatorname{diam}(T)$ and $f$ is a minimal dominating broadcast on
$T-b$. But then $T-b$ satisfies Theorem \ref{Thm_Cater}$(i)$ -- $(iii)$ and
$\Gamma_{b}(T-b)>\operatorname{diam}(T-b)$, contradicting the choice of $T$.
Hence assume $b$ is not a leaf adjacent to $v_{1}$.

\begin{figure}[ptb]%
\centering
\includegraphics[
height=0.8475in,
width=5.361in
]%
{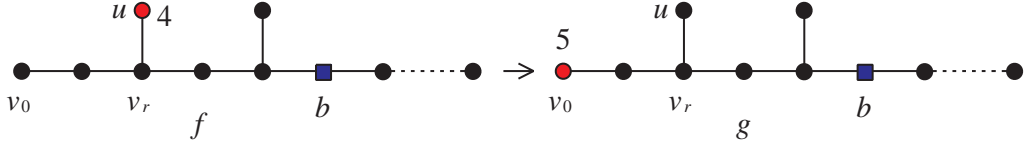}%
\caption{A step in the proof of Lemma \ref{Lem_non-diam1}.}%
\label{Fig_Cater2}%
\end{figure}

Let $r\geq1$ be the largest index such that $v_{r}$ lies on the $u-v_{0}$ path
in $T$. Possibly $v_{r}=u$, otherwise $u$ is a leaf adjacent to $v_{r}$. Since
$v_{0}$ is a peripheral vertex, $u$ broadcasts to all vertices of $T_{i}$ for
each $i=0,\ldots,r$, and each vertex $x$ in each such $T_{i}$ is overdominated
by $u$. Therefore $b\in V(T_{t})$ for some $t>r$. In addition, if $b$ lies on
$P$, then $b$ is not a stem, otherwise the leaves adjacent to $b$ are not
$f$-dominated. Therefore $u$ also broadcasts to each vertex of each $T_{i}$
for $r\leq i\leq t$. See Figure \ref{Fig_Cater2}.\emph{ }But then the
broadcast $g$ defined by $g(v_{0})=f(u)-d(u,v_{r})+d(v_{0},v_{r}),\ g(u)=0$
and $g(x)=f(x)$ otherwise is also a dominating broadcast such that
$b\in\operatorname{PB}_{g}(v_{0})$ and $\operatorname{PB}_{g}%
(x)=\operatorname{PB}_{f}(x)$ for all $x\in V_{g}^{+}-\{v_{0}\}$, that is, $g$
is a minimal dominating broadcast. Now $\sigma(g)\geq\sigma(f)-d(u,v_{r}%
)+d(v_{0},v_{r})\geq\sigma(f)-1+1=\sigma(f)$. Hence $\sigma(g)=\sigma(f)$ if
and only if $r=1$ and $u$ is a leaf adjacent to $v_{1}$. In this case,
$T-v_{0}$ also satisfies $(i)$ -- $(iii),$ $\operatorname{diam}(T-v_{0}%
)=\operatorname{diam}(T)$, and $f$ is also a minimal dominating broadcast on
$T-v_{0}$, contradicting the choice of $T$. Hence $\sigma(g)>\sigma(f)$ and we
again have a contradiction, because $\sigma(f)=\Gamma_{b}(T)$ and no minimal
dominating broadcast has cost greater than $\Gamma_{b}(T)$. This proves the
lemma for $v_{0}$. The result for $v_{d}$ follows by symmetry.~$\blacksquare$%

\begin{lemma}
\label{Lem_non-diam2}Let $T$ be a smallest non-diametrical caterpillar that
satisfies Theorem $\ref{Thm_Cater}(i)$ -- $(iii)$ and $f$ be a minimal
dominating broadcast on $T$ such that $\sigma(f)>\operatorname{diam}(T)$. Then
each leaf $w\notin\{v_{0},v_{d}\}$ of $T$ is either a broadcast vertex or
$\operatorname{PB}_{f}(u)=\{w\}$ for some $u\in V_{f}^{+}$.
\end{lemma}

\noindent\textbf{Proof.\hspace{0.1in}}Suppose the conclusion is false and
$w\notin\{v_{0},v_{d}\}$ is a leaf of $T$ that is neither a broadcast vertex
nor the only vertex in the private boundary of some $u\in V_{f}^{+}$. Then
$T-w$ is a tree with diameter $d$ that satisfies $(i)$ -- $(iii)$, and $f$ is
a minimal dominating broadcast on $T-w$ as well, contrary to the choice of
$T$.~$\blacksquare$\medskip

We now return to Lemmas \ref{Lem_non-diam3} and \ref{Lem_K2}.\medskip

\noindent\textbf{Lemma \ref{Lem_non-diam3}\hspace{0.1in}}\emph{If }$T$\emph{
is a smallest non-diametrical caterpillar that satisfies Theorem
}$\ref{Thm_Cater}(i)$\emph{ -- }$(iii)$\emph{, then no vertex of }$T$\emph{ is
a strong stem.\medskip}

\noindent\textbf{Proof.\hspace{0.1in}}Suppose, to the contrary, that some
vertex $v$ of $T$ is a strong stem. Then $v=v_{i}$ for some $i$, since $T$ is
a caterpillar. Say $v_{i}$ is adjacent to the leaves $\ell$ and $\ell^{\prime
}$. Let $f$ be a minimal dominating broadcast on $T$ such that $\sigma
(f)>\operatorname{diam}(T)$. By Lemmas \ref{Lem_non-diam1} and
\ref{Lem_non-diam2} we may assume that each leaf of $T$ is either a broadcast
vertex, or the only vertex in the $f$-private boundary of some vertex in
$V_{f}^{+}$. Let $u$ be the vertex that broadcasts to $\ell$.

Suppose $u\neq\ell$. Then $\operatorname{PB}_{f}(u)=\{\ell\}$. If $u\neq
\ell^{\prime}$, then $d(u,\ell)=d(u,\ell^{\prime})$ and we also have
$\ell^{\prime}\in\operatorname{PB}_{f}(u)$, contrary to Lemma
\ref{Lem_non-diam2}. Hence $u=\ell^{\prime},\ f(\ell^{\prime})=2$ and
$\operatorname{PB}_{f}(\ell^{\prime})=\{\ell\}$. Let $H_{1}$ and $H_{2}$ be
the subtrees of $T-v_{i}$ that contain $v_{0}$ and $v_{d}$, respectively. If
$i\in\{1,d-1\}$, assume without loss of generality that $i=d-1$ and ignore
$H_{2}$. Since $\operatorname{diam}(T)\geq4$, $H_{1}$ is nontrivial. By
Theorem \ref{Thm_Cater}$(ii)$, $v_{i-1}$ is not a stem of $T$, hence
$\operatorname{diam}(H_{1})=i-1$. Since $\ell^{\prime}$ broadcasts to
$v_{i-1}$ and $\operatorname{PB}_{f}(\ell^{\prime})=\{\ell\}$, $v_{i-1}$ also
hears $f$ from some vertex $w\in V_{f}^{+}-\{\ell^{\prime}\}$. Since $\ell
\in\operatorname{PB}_{f}(\ell^{\prime})$, $w\notin\{v_{i},\ell\}\cup V(H_{2}%
)$, hence $w\in V(H_{1})$. By Proposition \ref{Prop_Disjoint} applied to
$w,\ell^{\prime}\in V_{f}^{+}$, $\operatorname{PB}_{f}(w)\subseteq V(H_{1})$.
Therefore $f\upharpoonleft H_{1}$ is a minimal dominating broadcast on $H_{1}$.

\begin{itemize}
\item If $v_{i-2}$ is not a stem of $T$, then either $v_{i-2}$ is adjacent to
only one leaf in $H_{1}$, namely $v_{i-1}$, in which case $H_{1}$ satisfies
Theorem \ref{Thm_Cater}$(i)$ -- $(iii)$, or $v_{i-2}$ is adjacent to the two
leaves $v_{i-1}$ and $v_{0}$ in $H_{1}$, in which case $H_{1}\cong P_{3}$.

\item On the other hand, if $v_{i-2}$ is a stem of $T$, then by Theorem
\ref{Thm_Cater}$(iii)$ and the fact that $v_{i}$ is adjacent to two leaves,
$v_{i-2}$ is adjacent to exactly one leaf in $T$, so that it is adjacent to
two leaves in $H_{1}$. If $v_{i-2}$ is the only strong stem of $H_{1}$, then
$H_{1}$ satisfies Theorem \ref{Thm_Cater}$(i)$ -- $(iii)$. Hence suppose that
for some $i^{\prime}<i-2$, $v_{i^{\prime}}$ is a strong stem (of $H_{1}$ and
of $T$). Since $(iii)$ holds for $T$, and $\deg_{T}(v_{i-2}),\deg_{T}%
(v_{i})>2$, there exists an index $k,\ i^{\prime}<k<i-2$, such that $\deg
_{T}(v_{k})=\deg_{T}(v_{k+1})=2$. Therefore $H_{1}$ satisfies Theorem
\ref{Thm_Cater}$(i)$ -- $(iii)$ in this case as well.
\end{itemize}

By the choice of $T$, $\Gamma_{b}(H_{1})=\operatorname{diam}(H_{1})=i-1$ in
all cases. Since $f\upharpoonleft H_{1}$ is a minimal dominating broadcast on
$H_{1}$, $\sigma(f\upharpoonleft H_{1})\leq i-1$. Similarly, if $1<i<d$,
$H_{2}$ satisfies Theorem \ref{Thm_Cater}$(i)$ -- $(iii)$ and $f\upharpoonleft
H_{2}$ is a minimal dominating function of $H_{2}$ such that $\sigma
(f\upharpoonleft H_{2})\leq\operatorname{diam}(H_{2})=d-i-1$. But then
$\sigma(f)=\sigma(f\upharpoonleft H_{1})+\sigma(f\upharpoonleft H_{2})+2\leq
d$ (or $\sigma(f)=\sigma(f\upharpoonleft H_{1})+2\leq d-2+2=d$, if $i=d-1$),
which is a contradiction because $T$ is non-diametrical.

Hence we may assume that $u=\ell$; that is, $\ell$ is a broadcast vertex. If
$\ell$ broadcasts to $\ell^{\prime}$, we get a contradiction as above. Hence
$f(\ell)=1=f(\ell^{\prime})$ (since no other vertex can broadcast to
$\ell^{\prime}$ without broadcasting to $\ell$). Then $v_{i}\notin%
\operatorname{PB}_{f}(x)$ for each $x\in V_{f}^{+}$. We may now define $H_{1}$
and $H_{2}$ as above and proceed as before to obtain a
contradiction.~$\blacksquare$

\medskip

\noindent\textbf{Lemma \ref{Lem_K2}\hspace{0.1in}}\emph{If }$T$\emph{ is a
smallest non-diametrical caterpillar that satisfies Theorem }$\ref{Thm_Cater}%
(i)$\emph{ -- }$(iii)$\emph{, then no vertex }$v_{i},\ i\in\{2,\ldots
,d-2\}$\emph{, is adjacent to a leaf.\medskip}

\noindent\textbf{Proof.\hspace{0.1in}}Suppose, to the contrary, that some
$v_{i},\ i\in\{2,\ldots,d-2\}$\emph{,} is adjacent to a leaf and let $k$ be
the largest index in $\{2,\ldots,d-2\}$ such that $v_{k}$ is a stem. By Lemma
\ref{Lem_non-diam3} we may assume that $T$ has no strong stems. By Theorem
\ref{Thm_Cater}$(ii)$, $\deg_{T}(v_{k-1})=\deg_{T}(v_{k+1})=2$. Let $\ell$ be
the leaf adjacent to $v_{k}$ and let $f$ be a minimal dominating broadcast on
$T$ such that $\sigma(f)>\operatorname{diam}(T)$. By Lemmas
\ref{Lem_non-diam1} and \ref{Lem_non-diam2} we may assume that each of $\ell$
and $v_{d}$ is either a broadcast vertex or the only vertex in the $f$-private
boundary of some vertex in $V_{f}^{+}$. We consider several cases. In each
case we delete an edge to obtain subtrees of $T$, each of which contains at
most one strong stem. Since $T$ satisfies Theorem \ref{Thm_Cater}$(i)$ --
$(iii)$, so do the subtrees. By the choice of $T$, each subtree thus obtained
is diametrical. We omit these details in the cases for the sake of brevity.

\medskip

\noindent\textbf{Case 1\hspace{0.1in}}$\ell$ belongs to a private boundary and
$v_{d}\in V_{f}^{+}$. Then either $\ell\in V_{f}^{+}$ and $\ell\in
\operatorname{PB}_{f}(\ell)$, or $\operatorname{PB}_{f}(u)=\{\ell\}$ for a
vertex $u\neq\ell$.\medskip

\noindent\textbf{Case 1(a)\hspace{0.1in}}$\{\ell\}=\operatorname{PB}_{f}%
(v_{d})$. Then $f(v_{d})=d-k+1$ and $v_{d}$ broadcasts to $v_{k-1}$. Hence
$v_{k-1}$ does not belong to the private boundary of any vertex in $V_{f}^{+}%
$. Therefore $v_{k-1}$ also hears $f$ from a vertex in $V_{f}^{+}-\{v_{d}\}$.
Also, $\{v_{k},\ldots,v_{d-1}\}\cap V_{f}^{+}=\varnothing$. Let $T^{\prime}$
be the subtree of $T-v_{k-1}v_{k}$ that contains $v_{0}$. For each vertex
$u\in V_{f}^{+}\cap V(T^{\prime})$, Proposition \ref{Prop_Disjoint} applied to
$u$ and $v_{d}$ implies that $\operatorname{PB}_{f}(u)\subseteq V(T^{\prime}%
)$. Therefore $f\upharpoonleft T^{\prime}$ is a minimal dominating broadcast
on $T^{\prime}$. By the choice of $T$, $\sigma(f\upharpoonleft T^{\prime}%
)\leq\operatorname{diam}(T^{\prime})=k-1$. But now $\sigma(f)=\sigma
(f\upharpoonleft T^{\prime})+f(v_{d})\leq k-1+d-k+1=d$, a contradiction.

\medskip

\noindent\textbf{Case 1(b)\hspace{0.1in}}$\ell\in\operatorname{PB}%
_{f}(u),\ u\neq v_{d}$ (possibly $u=\ell$). Then $u$ broadcasts to $v_{k}$,
hence $v_{k}\notin\operatorname{PB}_{f}(v_{d})$. By Proposition
\ref{Prop_Disjoint} and the choice of $k$ as the largest index such that
$v_{k}\neq v_{d-1}$ is a stem, there exists an index $j>k$ such that $v_{j}%
\in\operatorname{PB}_{f}(v_{d})$ (and thus $f(v_{d})=d-j$). Evidently, then,
the edge $v_{j-1}v_{j}$ does not hear $f$ from any vertex in $V_{f}^{+}$. Let
$T^{\prime}$ be the subtree of $T-v_{j-1}v_{j}$ that contains $v_{0}$. As in
Case 1(a) we see that $f\upharpoonleft T^{\prime}$ is a minimal dominating
broadcast on $T^{\prime}$.

If $j=k+1$, then $u$ broadcasts to $v_{k}$ and $\ell$ but not to $v_{k+1}$.
(This is only possible if $u=\ell$ and $f(\ell)=1$.) In this case,
$\operatorname{diam}(T^{\prime})=k+1$ and $f(v_{d})=d-j=d-k-1$.

If $j>k+1$, i.e., $j-1\geq k+1$, then $\operatorname{diam}(T^{\prime})=j-1$.
In either case we obtain a contradiction as before as in Case 1(a).

\medskip

\noindent\textbf{Case 2\hspace{0.1in}}$\ell\in\operatorname{PB}_{f}(u)$ and
$v_{d}\in\operatorname{PB}_{f}(w)$. By Lemma \ref{Lem_non-diam1}, $u\neq w$.
\medskip

\noindent\textbf{Case 2(a)\hspace{0.1in}}$v_{d}\in\operatorname{PB}_{f}%
(v_{d})$. Then $f(v_{d})=1$. If $\operatorname{PB}_{f}(v_{d})=\{v_{d-1}%
,v_{d}\}$, delete the edge $v_{d-2}v_{d-1}$ and proceed as in Case 1(b) to get
a contradiction. Thus, assume $\operatorname{PB}_{f}(v_{d})=\{v_{d}\}$. Then
$v_{d-1}$ hears $f$ from some other vertex as well, hence $f\upharpoonleft
(T-v_{d})$ is a minimal dominating broadcast on $T-v_{d}$. By the choice of
$T$, $\sigma(f\upharpoonleft(T-v_{d}))\leq d-1$ and so $\sigma(f)\leq d$, a
contradiction. \medskip

\noindent\textbf{Case 2(b)\hspace{0.1in}}$\{v_{d}\}=\operatorname{PB}_{f}(w)$
for some $w\neq v_{d}$. Since $w$ does not broadcast to $\ell$, Proposition
\ref{Prop_Disjoint} and the choice of $k$ imply that $w=v_{i}$ for some $i\geq
k+1$. Since $v_{d}\in\operatorname{PB}_{f}(v_{i})$, $f(v_{i})=d-i$. Let
$r=\min\{i,\min\{j:v_{j}\in\operatorname{PN}_{f}(v_{i})\}$. Since $u$
broadcasts to $v_{k}$ and $i\geq k+1$, $r\geq k+1$. Let $T^{\prime}$ and
$T^{\prime\prime}$ be the subtrees of $T-v_{r-1}v_{r}$ that contain $v_{0}$
and $v_{d}$, respectively. Then $\operatorname{diam}(T^{\prime\prime})=d-r$
and $V_{f}^{+}\cap V(T^{\prime\prime})=\{v_{i}\}$. Since $r\leq i$,
$f(v_{i})=d-i\leq d-r$.

If $r=i$, then $v_{r},\ldots,v_{d}$ is a path from $v_{r}$ to $v_{d}%
\in\operatorname{PB}_{f}(v_{r})$. Otherwise, $r<i$ and, by definition of $r$,
$\{v_{r},\ldots,v_{i},\ldots,v_{d}\}\subseteq\operatorname{PN}_{f}(v_{i})$. In
either case, Proposition \ref{Prop_Disjoint} again implies that
$\operatorname{PB}_{f}(x)\subseteq V(T^{\prime})$ for each $x\in V_{f}%
^{+}-\{v_{i}\}$. Hence $f\upharpoonleft T^{\prime}$ is a minimal dominating
broadcast on $T^{\prime}$, so that by the choice of $T$, $\sigma
(f\upharpoonleft T^{\prime})\leq\operatorname{diam}(T^{\prime})$. If $r=k+1$,
then $\operatorname{diam}(T^{\prime})=r$, while if $r>k+1$, then
$\operatorname{diam}(T^{\prime})=r-1$. In either case $\sigma(f)\leq
r+f(v_{i})\leq d$, a contradiction.

\medskip

\noindent\textbf{Case 3\hspace{0.1in}}$\ell$ is a broadcast vertex and
$\operatorname{PB}_{f}(u)=\{v_{d}\}$ for some vertex $u\neq v_{d}$. \medskip

\noindent\textbf{Case 3(a)\hspace{0.1in}}$\operatorname{PB}_{f}(\ell
)=\{v_{d}\}$. Then $f(\ell)=d-k+1\geq3$. Let $P$ be the $\ell-v_{d}$ path in
$T$ and let $w\in V_{f}^{+}-\{\ell\}$. Then $P\cong P_{f(\ell)+1}$. By
Proposition \ref{Prop_Disjoint}, $w\in V(T_{i})$ for some $i\leq k-1$. Also,
$\operatorname{PB}_{f}(w)\cap V(P)=\varnothing$. Thus, if $T^{\prime}$ is the
subtree of $T-v_{k-1}v_{k}$ that contains $v_{0}$, then $\operatorname{diam}%
(T^{\prime})=k-1$ and $f\upharpoonleft T^{\prime}$ is a minimal dominating
broadcast on $T^{\prime}$, which is a diametrical tree. Now $\sigma
(f)=\sigma(f\upharpoonleft T^{\prime})+f(\ell)\leq k-1+d-k+1=d$, a
contradiction. \medskip

\noindent\textbf{Case 3(b)\hspace{0.1in}}$w\in\operatorname{PB}_{f}(\ell)$ and
$\operatorname{PB}_{f}(u)=\{v_{d}\}$, where $u\notin\{\ell,v_{d}\}$ and $w\neq
v_{d}$. By Proposition \ref{Prop_Disjoint}, $u=v_{i}$ for some $i\geq k+1$. We
now proceed as in Case 2(b) to obtain a contradiction.

\medskip

\noindent\textbf{Case 4\hspace{0.1in}}$\ell$ and $v_{d}$ are both broadcast
vertices. If $f(v_{d})=1$, then $v_{d}\in\operatorname{PB}_{f}(v_{d})$. This
is Case 2(a), hence assume $f(v_{d})\geq2$. Then $\operatorname{PB}_{f}%
(v_{d})=\{v_{i}\}$ for some $i$ such that $k+1\leq i\leq d-2$. Evidently,
then, the edge $e=v_{i-1}v_{i}$ does not hear $f$ from any vertex. By deleting
$e$ we proceed as before to obtain a contradiction.\medskip

Since Cases 1 -- 4 and their subcases cover all possibilities for $\ell$ and
$v_{d}$, the lemma follows.~$\blacksquare$

\bigskip

This concludes the proofs of Lemmas \ref{Lem1-1} -- \ref{Lem_K2}, hence the
proof of Theorem \ref{Thm_Cater} is complete.

\section{Open Problems\label{Sec_Problems}}

A characterization of diametrical caterpillars is presented in Theorem
\ref{Thm_Cater}. In general, diametrical trees can have paths of length one or
two, but not longer paths, that are internally disjoint from a diametrical path.

\begin{problem}
Characterize diametrical trees that contain at least one path of length two
internally disjoint from a diametrical path.
\end{problem}

\begin{problem}
Characterize trees $T$ with $(i)$ $\Gamma_{b}(T)=\alpha(T),\ (ii)\ \Gamma
_{b}(T)=\Gamma(T)$.
\end{problem}

\begin{problem}
Study other classes of graphs $G$ such that $\Gamma_{b}(G)=\operatorname{diam}%
(G)$.
\end{problem}

\begin{problem}
\emph{\cite{MR}}\hspace{0.1in}Determine the maximum ratio $\Gamma
_{b}(G)/\Gamma(G)$ for $(i)$ general graphs, $(ii)$ trees.
\end{problem}

The stars $K_{1,n}$ satisfy $\operatorname{diam}(K_{1,n})=2$ and $\Gamma
_{b}(K_{1,n})=n$, hence the ratio $\Gamma_{b}(G)/\operatorname{diam}(G)$ is unbounded.

The proof of Lemma \ref{Lem1-1} suggests the following problem.

\begin{problem}
\label{ref}If $G$ and $H$ are graphs and $G$ is an isometric subgraph of $H$,
is it true that $\Gamma_{b}(G)\leq\Gamma_{b}(H)$?
\end{problem}

\subsection*{Acknowledgements}

The authors are indebted to the referees for several corrections and
improvements to the paper. In particular, Problem \ref{ref} was suggested by
one of them.

\bigskip


\begin{thebibliography}{99}                                                                                               %


\bibitem {Ahmadi}D.~Ahmadi, G.~H.~Fricke, C.~Schroeder, S.~T.~Hedetniemi,
R.~C.~Laskar, Broadcast irredundance in graphs. \emph{Congr.~Numer.}%
~\textbf{224} (2015), 17--31.

\bibitem {BF}I. Bouchemakh and N. Fergani, On the upper broadcast domination
number, \emph{Ars Combin}. \textbf{130} (2017), 151-161.

\bibitem {BouchS}I.~Bouchemakh, R.~Sahbi, On a conjecture of Erwin,
\emph{Stud. Inform. Univ.} \textbf{9}(2) (2011), 144--151.

\bibitem {BZ}I.~Bouchemakh, M.~Zemir, On the broadcast independence number of
grid graph, \emph{Graphs Combin}. \textbf{30} (2014),83--100.

\bibitem {BMT}R.~C.~Brewster, C.~M.~Mynhardt, L.~Teshima, New bounds for the
broadcast domination number of a graph, \emph{Central European J. Math.}
\textbf{11}(7) (2013), 1334--1343.

\bibitem {BS}B.~Bre\v{s}ar, S.~\v{S}pacapan, Broadcast domination of products
of graphs, \emph{Ars Combin}. \textbf{92} (2009), 303--320.

\bibitem {CLZ}G.~Chartrand, L.~Lesniak, P.~Zhang, \emph{Graphs and Digraphs}
(Sixth Edition), Chapman \&\ Hall, 2015.

\bibitem {CHM}E.~J.~Cockayne, S.~Herke, C.~M.~Mynhardt, Broadcasts and
domination in trees, \emph{Discrete Math.} \textbf{311} (2011), 1235--1246.

\bibitem {DDH}J.~Dabney, B.~C.~Dean, S.~T.~Hedetniemi, A linear-time algorithm
for broadcast domination in a tree, \emph{Networks }\textbf{53} (2009) 160--169.

\bibitem {Dunbar}J.~Dunbar, D.~Erwin, T.~Haynes, S.~M.~Hedetniemi,
S.~T.~Hedetniemi, Broadcasts in graphs, \emph{Discrete Applied Math.}
\textbf{154} (2006), 59-75.

\bibitem {Ethesis}D.~Erwin, \emph{Cost domination in graphs}. Doctoral
dissertation, Western Michigan University, 2001.

\bibitem {Epaper}D.~Erwin, Dominating broadcasts in graphs, \emph{Bulletin of
the ICA} \textbf{42} (2004), 89--105.

\bibitem {HHS}T.~W.~Haynes, S.~T.~Hedetniemi and P.~J.~Slater,
\emph{Fundamentals of Domination in Graphs}. Marcel Dekker, New York, 1998.

\bibitem {HL}P.~Heggernes, D.~Lokshtanov, Optimal broadcast domination in
polynomial time, \emph{Discrete Math.} \textbf{36} (2006), 3267-3280.

\bibitem {HS}P.~Heggernes, S.~H.~S\ae ther, Broadcast domination on block
graphs in linear time. Computer science -- theory and applications, 172--183,
\emph{Lecture Notes in Comput. Sci.} \textbf{7353}, Springer, Heidelberg, 2012.

\bibitem {Herke}S.~Herke, \emph{Dominating broadcasts in graphs}, Master's
thesis, University of Victoria, 2009. http://hdl.handle.net/1828/1479

\bibitem {HM}S.~Herke, C.~M.~Mynhardt, Radial Trees. \emph{Discrete Math.
}\textbf{309} (2009), 5950--5962.

\bibitem {RadKhos}N.~Jafari Rad, F.~Khosravi, Limited dominating broadcast in
graphs, \emph{Discrete Math. Algorithms Appl}. \textbf{5} (2013) [9 pages].
DOI: 10.1142/S1793830913500250

\bibitem {Scott}S.~Lunney, Trees with equal broadcast and domination numbers,
\emph{Master's thesis, University of Victoria}, 2011. http://hdl.handle.net/1828/3746

\bibitem {SM}S.~Lunney and C.~M.~Mynhardt, More trees with equal broadcast and
domination numbers, \emph{Australas.~J.~Combin.} \textbf{61} (2015), 251--272.

\bibitem {MR}C.~M.~Mynhardt, A.~Roux, Dominating and irredundant broadcasts in
graphs, submitted. ArXiv link http://arxiv.org/abs/1608.00052.

\bibitem {MT}C.~M.~Mynhardt, L.~Teshima, Broadcasts and multipackings in
trees, \emph{Utilitas Math.}, to appear.

\bibitem {MW}C.~M.~Mynhardt, J.~S.~Wodlinger, A class of trees with equal
broadcast and domination numbers. \emph{Australasian J. Math. }\textbf{56}
(2013), 3--22.

\bibitem {MW2}C.~M.~Mynhardt, J.~S.~Wodlinger, Uniquely radial trees, \emph{J.
Combin. Math. Combin. Comput.}, to appear. Accepted 19 May 2014.

\bibitem {SS}S.~M.~Seager, Dominating broadcasts of caterpillars, \emph{Ars
Combin.} \textbf{88} (2008), 307--319.
\end{thebibliography}
\end{document}